\documentclass[a4paper,12pt]{amsart}
\usepackage{graphicx}
\usepackage{amscd}
\usepackage{amssymb}
\newtheorem{thm}{Theorem}[section]
\newtheorem{cor}[thm]{Corollary}
\newtheorem{lem}[thm]{Lemma}
\newtheorem{prop}[thm]{Proposition}
\theoremstyle{defn}
\newtheorem{defn}[thm]{Definition}
\theoremstyle{rem}

\newtheorem{nota}[thm]{Notation}
\numberwithin{equation}{section}

\newcommand{\eps}{\varepsilon}

\newcommand{\G}{\mathcal{G}}
\newcommand{\M}{\mathcal{M}}

\newcommand{\EE}{\mathcal E_\pi}
\newcommand{\Euv}{\mathcal E_{u,v}^\pi}
\newcommand{\R}{\mathcal{R}ep(\mathcal G)}

\newcommand{\HH}{\mathcal H_\pi}
\newcommand{\Hu}{\mathcal H_u^\pi}
\newcommand{\Hv}{\mathcal H_v^\pi}
\newcommand{\PP}{(\pi,\mathcal H_\pi , \mu_\pi)}
\newcommand{\repn}{representation\,}
\newcommand{\repns}{representations\,}

\begin{document}

\title[duality for groupoids]{Tannaka-Krein duality for compact groupoids I,
representation theory}%
\author{Massoud Amini}%
\address{
Department of Mathematics, Tarbiat Modarres University, P.O.Box
14115-175, Tehran , Iran , mamini@modares.ac.ir \newline
Department of Mathematics and Statistics, University of
Saskatchewan, 106 Wiggins Road, Saskatoon, S7N 5E6 ,
mamini@math.usask.ca}

\thanks{The author was visiting the University of Saskatchewan during the
preparation of this work. He would like to thank University of
Saskatchewan and Professor Mahmood Khoshkam for their hospitality
and support. Also I would like to thank Professor Alan Paterson
who kindly pointed out an error in the preliminary version of this
paper and Professor Paul Muhly for sharing some of his unpublished
results with the author}%


\begin{abstract}
In a series of papers, we have shown that from the \repn theory of
a compact groupoid one can reconstruct the groupoid using the
procedure similar to the Tannaka-Krein duality for compact groups.
In this part we study continuous \repns of compact groupoids. We
show that irreducible \repns have finite dimensional fibres. We
prove the Schur's lemma, Gelfand-Raikov theorem and Peter-Weyl
theorem for compact groupoids.
\end{abstract}
\maketitle
\section{introduction}
The duality theory for Abelian locally compact groups [R] was
introduced by Pontyagin in 1934. Since then, many attempts are
done to generalize this duality theory for non Abelian locally
compact groups (see [ES] for a brief history and references). The
dual group in pontryagin theory is the group of characters with
pointwise multiplication. This group is not large enough to
recover the original group in the non Abelian case (there are
examples of non Abelian groups with a trivial group of
characters). The natural candidate in the non Abelian setting is
the set of (equivalence classes) of irreducible unitary
representations. This object is not a group, but it is shown to
recover the original group (at least in the compact case). One of
the successful methods of recovery was introduced by Tannaka in
1939 (see [JS] for a very clear exposition from a Category point
of view). Tannaka showed us how to recover a compact group from
the set (category) of its \repns . This is loosely called the
Tannaka duality (it is not a duality in the technical sense, it is
indeed an equivalence of categories).

Topological groupoids are natural generalizations of topological
groups [Re]. These are very rich structures and arise in a vast
spectrum of applications [P]. The \repn theory of groupoids is
more involved. In the group case, we represent group elements as
unitary operators on a (usually infinite dimensional) Hilbert
space. For groupoids we need a bundle of Hilbert spaces on the
unit space of the groupoid and each element is represented by a
bundle operator which is not unitary in general (the closest thing
you can imagine is a partial isometry). The concepts like the
group algebra, group $C^*$-algebra and regular \repn could
naturally be defined in this setting [Re], [P].

The main objective of this and the forthcoming papers [A1], [A2]
is to generalize the Tannaka duality to compact groupoids. The
next section of this paper includes the basic \repn theory of
compact groupoids. Some of the materials in this section are new.
The main result of this section asserts that irreducible
(continuous) representations of compact groupoid have finite
dimensional fibres . In contrast with compact groups, these \repns
are not necessarily finite dimensional (the integral of a finite
dimensional bundle which is irreducible in our sense is
irreducible in the usual sense , but could be infinite
dimensional). In the third section, we generalize the main results
of the Harmonic Analysis on compact groups to compact groupoids.
These include Schur's orthogonality relations , Schur's lemma,
Peter-Weyl theorem and Gelfand-Raikov theorem. In [A1] we
introduce the Fourier and Fourier-Plancherel transforms for
compact groupoids and study their properties. The main result of
these series of papers is the Tannaka-Krein theorem for compact
groupoids which is proved in [A2]. The proofs of most of the
results in sections two and three of the current follow closely
the analogous results for compact groups (see for instance [F]).
We give a detailed proof only when the modifications are
substantial.

\section{Irreducible representations of compact groupoids}

In this section we review  the \repn theory of locally compact
groupoids. Then we restrict ourselves to the compact case, in
which we prove analogues of some of the classical results for
\repns of compact groups. In the beginning of this section we
assume that $\G$ is a locally compact (not necessarily Hausdorff)
groupoid and $X=\G^{(0)}$.

\begin{defn} A \repn of $\G$ is a triple $\PP$ , where
$\HH=\{\mathcal H_u^\pi\}$ is a bundle of Hilbert spaces,
$\mu_\pi$ is a quasi-invariant measure on $X$ (with associated
measures $\nu, \nu^{-1}, \nu^2$, and $\nu_0$) such that

$(i) \, \pi(x)\in\mathcal B(\mathcal H_{s(x)}^\pi,\mathcal
H_{r(x)}^\pi)\quad (x\in \G)$,

$(ii)\, \pi(u)=id_u:\mathcal H_u^\pi\to\mathcal H_u^\pi\quad (u\in
X)$,

$(iii) \, \pi(xy)=\pi(x)\pi(y)$ for $\nu^2$-a.e.
$(x,y)\in\G^{(2)}$,

$(iv) \, \pi(x^{-1})=\pi(x)^{-1}$ for $\nu$-a.e. $x\in\G$,

$(v)\, x\mapsto <\pi(x)\xi(s(x)),\eta(r(x))>$ is $\nu$-measurable
on $\G$ for each $\xi,\eta\in L^2(\G^{(0)},\mathcal
H_\pi,\mu_\pi)$.
\end{defn}

For our purposes, we need more restricted version of \repns .

\begin{defn} A continuous \repn of $\G$ is a pair $(\pi,\mathcal H_\pi)$ , where
$\HH=\{\mathcal H_u^\pi\}$ is a bundle of Hilbert spaces, such
that $(i)-(iv)$ above hold everywhere (instead of almost
everywhere) and moreover the maps in $(v)$ are continuous.
\end{defn}

It is easy to check that

\begin{lem} If $\pi$ is a continuous \repn of $\G$, then for each
$\xi,\eta\in\HH$ the map $x\mapsto<\pi(x)\xi_{s(x)},\eta_{r(x)}>$ is continuous.\qed
\end{lem}

\begin{defn} Two \repns $(\pi_1,\mathcal H_{\pi_1} , \mu_{\pi_1})$ and
$(\pi_2,\mathcal H_{\pi_2} , \mu_{\pi_2})$ are called unitarily
equivalent if measures $\mu_{\pi_1})$ and $\mu_{\pi_2})$ are
equivalent (each is absolutely continuous with respect to the
other) and there is a bundle $\{U_u\}_{u\in X}$ with $U_u\in \mathcal B(\mathcal
H_u^{\pi_1},\mathcal H_u^{\pi_2})$ a unitary operator such that
for $\nu_1$-a.e. $x\in\G$ the following diagram commutes

\begin{equation*}
\begin{CD}
\mathcal H_{s(x)}^{\pi_1}@>{\pi_1(x)}>>\mathcal H_{r(x)}^{\pi_1}\\
@V{U_{s(x)}}VV    @VV{U_{r(x)}}V\\
\mathcal H_{s(x)}^{\pi_2}@>>{\pi_2(x)}>\mathcal H_{r(x)}^{\pi_2}
\end{CD}
\end{equation*}

The definition of equivalence of continuous \repns is similar,
except that we require the commutativity of the diagram for all
$x\in\G$. To show that $\pi_1$ and $\pi_2$ are unitarily equivalent we write $\pi_1\sim\pi_2$.
\end{defn}

\begin{nota}
We denote by $\mathcal{R}ep(\G)$ the category consisting of
(equivalence classes of) continuous \repns of $\G$ as objects and
intertwining operators as morphisms, namely $h=\{h_u\}_{u\in X}\in
Mor(\pi_1,\pi_2)$ if $h_u\in\mathcal B(H_u^{\pi_1},\mathcal
H_u^{\pi_2})$ make the following diagram commutative for each
$x\in\G$.

\begin{equation*}
\begin{CD}
\mathcal H_{s(x)}^{\pi_1}@>{\pi_1(x)}>>\mathcal H_{r(x)}^{\pi_1}\\
@V{h_{s(x)}}VV    @VV{h_{r(x)}}V\\
\mathcal H_{s(x)}^{\pi_2}@>>{\pi_2(x)}>\mathcal H_{r(x)}^{\pi_2}
\end{CD}
\end{equation*}

\end{nota}

With  this notation, two \repns $\pi_1$ and $\pi_2$ are unitarily
equivalent if there is a unitary operator bundle in
$Mor(\pi_1,\pi_2)$. In general, it is clear that
$Mor(\pi_1,\pi_2)$ is a vector space under pointwise operations on
operator bundles. Also it is routin to check that

\begin{lem}
For each $\pi\in\R$, $Mor(\pi,\pi)$ is a unital involutive algebra.\qed
\end{lem}

\begin{defn}
If  $\pi\in\R$ and $\M\subseteq \HH$ is a closed nonzero subbundle
(i.e. $\M_u$ is a nonzero closed subspace of $\Hu$, for each $u\in
X$) which is invariant under $\pi$ (i.e. $\pi(x)\M_{s(x)}\subseteq
\M_{r(x)}$, for each $x\in\G$) then $\pi^\M
(x)=\pi(x)\upharpoonleft_{\M_{s(x)}}$ defines a \repn of $\G$ on
$\M$ which is called a sub\repn of $\pi$. If $\pi$ admits a
nontrivial invariant subbundle $\M$ (nontrivial means $\M\neq 0$
or $\HH$) it is called  reducible, otherwise it is called
irreducible. We denote the set of equivalence classes of
irreducible \repns of $\G$ by $\hat\G$.
\end{defn}

\begin{defn}
If  $\{\pi_i\}_{i\in I}$ is a family of (continuous) \repns of
$\G$, their direct sum $\pi=\bigoplus_{i\in I} \pi_i$ on
$\HH=\mathcal H_{\pi_i}$ is defined by
$$\pi(x)(\xi_i)_{i\in I}=(\pi_i(x)(\xi_i))_{i\in I}
\quad (\xi_i\in\mathcal H_{s(x)}^{\pi_i},\, i\in I).$$
We  regard $\mathcal H_{\pi_i}$ as an invariant closed subbundle
of $\HH$ and $\pi_i$ as a sub\repn of $\pi$.
\end{defn}

\begin{lem}
If $\M$ is an  invariant subbundle of $\HH$, then so is $\M^\perp
:=\{\M_u^\perp\}_{u\in X}$. Moreover the later is closed (even if
the former is not!).
\end{lem}
{\bf Proof} For  each $x\in \G$, $\xi\in \M_{s(x)}^\perp$, and
$\eta\in\M_{r(x)}$ we have
$<\pi(x)\xi,\eta>=<\xi,\pi(x^{-1})\eta>=0$, that is
$\pi(x)\M_{s(x)}^\perp\subseteq \M_{r(x)}^\perp$, so $\M^\perp$ is
invariant. The last statement follows from the same fact for
Hilbert spaces.\qed

\begin{defn}
Let $\pi\in\R$  and $\xi\in\HH$. The subbundle $\M_\xi$ whose leaf
at $u\in X$ is the closed linear span of the set
$\{\pi(x)\xi_{s(x)}: x\in \G^u\}$ in $\Hu$ is called the cyclic
subbundle generated by $\xi$. This is a closed invariant subbundle
for $\pi$. We say that $\xi$ is a cyclic vector for $\pi$, if
$(\M_\xi)_u$ is dense in $\Hu$, for each $u\in X$. In this case,
$\pi$ is called a cyclic \repn.
\end{defn}

Next result follows  from above lemma and a standard argument
based on Zorn's lemma (see [F, 3.3]).

\begin{prop}
Each (continuous) \repn of $\G$ is the direct sum of cyclic
(continuous) \repns .\qed
\end{prop}

\begin{lem}
Let $\pi\in\R$, $\M$  be a closed subbundle of $\HH$, and
$P:\HH\to\M$ be the corresponding orthogonal projection. Then $\M$
is invariant under $\pi$ if and only if $P\in Mor(\pi,\pi)$.
\end{lem}
{\bf Proof}
If  $P\in Mor(\pi,\pi)$, $\xi\in\M$, then for each $x\in\G$
$$\pi(x)\xi_{s(x)}=\pi(x)P_{s(x)}\xi_{s(x)}=P_{r(x)}\pi(x)\xi_{s(x)}\in\M_{r(x)},$$
so $\M$ is invariant. Conversely if $\M$  is invariant then , by
above lemma, for each $\xi\in\M$, $\eta\in\M^\perp$, and $x\in\G$
we have
$$\pi(x)P_{s(x)}\eta_{s(x)}=0=P_{r(x)}\pi(x)\eta_{s(x)},$$
so $P\in Mor(\pi,\pi)$. \qed

\vspace{.3 cm} Next we need to recall a standard result from the
theory of measurable functional calculus.

\begin{lem} If $\mathcal H$ and $\mathcal K$ are Hilbert spaces and
$T_1\in\mathcal B(\mathcal H)$ and $T_2\in\mathcal B(\mathcal K)$ are normal operators,
and $S\in\mathcal B(\mathcal H,\mathcal K)$
satisfies $ST_1=T_2S$, then for each bounded Borel map
$f:\sigma(T_1)\cup\sigma(T_2)\to\mathbb R$, $Sf(T_1)=f(T_2)S$.\qed
\end{lem}

Now we are ready to prove Schur's lemma for compact groupoids.

\begin{thm} {\bf (Schur's lemma)}
A  (continuous) \repn $\pi$ of $\G$ is irreducible if and only if
$Mor(\pi,\pi)\simeq\mathbb C$. If $\pi_1,\pi_2$ are irreducible
(continuous) \repns of $\G$, then
\begin{equation}
Mor(\pi_1,\pi_2)\simeq
\begin{cases}
\,\mathbb C& \text{if} \,\, \pi_1\sim\pi_2,\\
\,0& \text{otherwise}
\end{cases}
\end{equation}
\end{thm}
{\bf Proof} If $\pi$ is reducible,  then by above proposition
$Mor(\pi,\pi)$ contains a nontrivial projection bundle. Conversely
if $T\in Mor(\pi,\pi)$ is not a multiple of identity, then by
lemma 2.5, $A=\frac{1}{2}(T+T^*)$ and $B=\frac{1}{2i}(T-T^*)$ are
in $Mor(\pi,\pi)$ and at least one of them, say $A$ is not a
multiple of identity. $A$ is self-adjoint, so by above lemma,
applied to $T_1=A_{s(x)}, T_2=A_{r(x)}$, $S=\pi(x)$, and
$f=\chi_E$, where $E$ is a Borel subset of $\mathbb R$, we get
$$\pi(x)\chi_E(A_{s(x)})=\chi_E(A_{r(x)})\pi(x)\quad (x\in\G),$$
so if we put $\chi_E(A)=\{\chi_E(A_u)\}_{u\in X}$, then
$\chi_E(A)\in Mor(\pi,\pi)$,  thereby $Mor(\pi,\pi)$ contains at
least one nontrivial projection bundle, and so again by Lemma 2.5,
$\pi$ is reducible.

Next for irreducible \repns  $\pi_1,\pi_2$ of $\G$, let $T\in
Mor(\pi_1,\pi_2)$, then clearly $T^*\in Mor(\pi_2,\pi_1)$ and so
$T^*T\in Mor(\pi_1,\pi_1)$ and $TT^*\in Mor(\pi_2,\pi_2)$. Hence
$TT^*$ and $T^*T$ are both multiples of identity. So if $T\neq 0$,
then a multiple of $T$ is unitary. Therefore
$Mor(\pi_1.\pi_2)=\{0\}$, precisely when $\pi_1$ and $\pi_2$ are
not equivalent. Now if $T_1,T_2\in Mor(\pi_1,\pi_2)$ and $T_2\neq
0$, then $T_2$ is a (nonzero) multiple of a unitary, so
$T_2^{-1}T_1=T_2^*T_1\in Mor(\pi_1)$ is a multiple of identity, so
$T_1$ is a multiple of $T_2$. \qed

\vspace{.3 cm} In the rest of this section we  assume that $\G$ is
compact. Note that in this case for each $u\in X$, the subsets
$\G_u$ and $\G^u$ of $\G$ are compact. In particular the isotropy
groups $\G_u^u$ of $\G$ are compact groups. We may assume that the
Haar system of $\G$ is normalized in a way that
$\lambda_u(\G_u)=\lambda^u(\G^u)=1$, for each $u\in X$.

\begin{lem}
Assume that $\G$ is compact.  If $\pi\in\R$ is a continuous
irreducible \repn of $\G$, $\xi\in\HH$ with $\|\xi_u\|=1$, $u\in
X$, and $T_u:\Hu\to\Hu$ is defined by
$$T_u\eta_u=\int <\eta_{r(x)},\pi(x)\xi_{s(x)}>\pi(x)\xi_{s(x)}
d\lambda^u(x)\quad (\eta\in\HH),$$
Then $T\geq 0$, $T\neq 0$, and $T\in Mor(\pi,\pi)\cap\mathcal K(\HH)$.
\end{lem}
{\bf Proof}
For $\eta\in\HH$ and $u\in X$
\begin{align*}
<T_u\eta_u,\eta_u>&=\int<\eta_{r(x)},\pi(x)\xi_{s(x)}>
<\pi(x)\xi_{s(x)},\eta_{u}>d\lambda^u(x)\\
&=\int |<\eta_u,\pi(x)\xi_{s(x)}>|^2 d\lambda^u(x)\geq 0.
\end{align*}
This shows that $T\geq 0$ and when $\eta=\xi$, it tells us that the real valued function
$$f(x)=|<\xi_{s(x)},\pi(x)\xi_{s(x)}>|\quad(x\in \G^u)$$
is  non negative. Also we have $f(u)=\|\xi_u\|^2=1$. By Lemma 2.3,
$f$ is continuous, so there is a Hausdorff open neighborhood $V$
of $u$ in $\G^u$ such that $|f|>\frac{1}{2}$ on $V$. Therefore
$$<T_u\xi_u,\xi_u>=\int_{\G^u} f(x)d\lambda^u(x)\geq\int_{V} f(x)d\lambda^u(x)
\geq\frac{1}{2}\lambda^u(V)>0,$$
so $T_u\neq 0$.

Next  fix $\eta\in\HH$ with $\|\eta_u\|=1, u\in X$. Fix $u\in X$.
By compactness of $\G^u$, the map $x\mapsto \pi(x)\xi_s(x)$ is
uniformly continuous on $\G^u$, so given $\eps>0$, there is a
partition $E_1,\dots,E_n$ of $\G^u$ into finitely many mutually
disjoint Borel subsets and elements $x_j\in E_j\, (1\leq j\leq n)$
such that
$$\|\pi(x)\xi_{s(x)}-\pi(x_j)\xi_{s(x_j)}\|<\frac{\eps}{2},
\quad  (x\in E_j, 1\leq j\leq n).$$
For each $x\in E_j$
\begin{align*}
\|<\eta_{u},\pi(x)\xi_{s(x)}>\pi(x)\xi_{s(x)}
&-<\eta_{u},\pi(x_j)\xi_{s(x_j)}>\pi(x_j)\xi_{s(x_j)}\|\\
&\leq\|<\eta_{u},(\pi(x)\xi_{s(x)}-\pi(x_j)\xi_{s(x_j)}>\pi(x)\xi_{s(x)}\|\\
&+\|<\eta_{u},\pi(x_j)\xi_{s(x_j)}>(\pi(x)\xi_{s(x)}-\pi(x_j)\xi_{s(x_j)}\|\\
&\leq \frac{\eps}{2}(2\|\eta_{u}\|)=\eps,
\end{align*}
so if we put
\begin{align*}
T_u^\eps\eta_u&=\sum_{j=1}^{n} \lambda^u(E_j)
<\eta_{u},\pi(x_j)\xi_{s(x_j)}>\pi(x_j)\xi_{s(x_j)}\\
&=\sum_{j=1}^{n} \int_{E_j} <\eta_{u},\pi(x_j)\xi_{s(x_j)}>\pi(x_j)\xi_{s(x_j)} d\lambda^u(x),
\end{align*}
then we have
\begin{align*}
\|T_u\eta_u-T_u^\eps\eta_u\|&=\|\sum_{j=1}^{n} \int_{E_j}
<\eta_{u},\pi(x)\xi_{s(x)}>\pi(x)\xi_{s(x)} \\
&-<\eta_{u},\pi(x_j)\xi_{s(x_j)}>\pi(x_j)\xi_{s(x_j)}>d\lambda^u(x)\|<\eps.
\end{align*}
Therefore  $\|T_u^\eps-T_u\|<\eps$, and as each $T_u^\eps$ is
clearly a finite rank operator, each $T_u$ is a compact operator.
Finally $T\in Mor(\pi,\pi)$, because for each $x\in \G$,
$\eta\in\HH$, and $u\in X$
\begin{align*}
\pi(x)T_{s(x)}\eta_s(x)&=\int
<\eta_{s(x)},\pi(y)\xi_{s(y)}>\pi(x)\pi(y)\xi_{s(y)} d\lambda^{s(x)}(y)\\
&=\int <\eta_{s(x)},\pi(x^{-1}y)\xi_{s(y)}>\pi(y)\xi_{s(y)}
d\lambda^{r(x)}(y)=T_{r(x)}\pi(x)\eta_{s(x)}.\qed
\end{align*}

\vspace{.3 cm} The above lemma enables  us to prove the finite
dimensionality of bundles of continuous irreducible \repns of
compact groupoids. The proof now goes exactly as in the group case
[F, 5.2]. we give the proof for the sake of completeness.

\begin{thm}
If $\G$ is compact,  then each  continuous irreducible \repn of
$\G$ has a finite dimensional bundle and each continuous \repn of
$G$ is a direct sum of continuous irreducible \repns .
\end{thm}
{\bf Proof} If $\pi$ is  irreducible and $T$ is as above, then
each $T_u$ is compact and a nonzero multiple of the identity on
$\Hu$, so $dim(\Hu)<\infty$, for each $u\in X$. Now if $\pi$ is
arbitrary, then each $T_u$ has a nonzero eigenvalue in
$\zeta_u\in\Hu$ [F, 1.52] whose eigenspace $\M_u$ is finite
dimensional. Let $\M=\{\M_u\}$, then $\M$ is invariant under $T$,
since $T\in Mor(\pi,\pi)$, hence $\pi$ has a finite dimensional
sub\repn . But an iductive argument based on Lemma 2.9 shows that
any finite dimensional \repn has an irreducible sub\repn , and so
does $\pi$ as well. Now by Zorn Lemma we can find a maximal family
of mutually orthogonal irreducible invariant subbundles of $\HH$,
whose direct sum has to be $\HH$ by Lemma 2.9 and maximality.\qed

\vspace{.3 cm} Note that  the decomposition into irreducible
\repns is not in general unique, but the decomposition into
subspaces corresponding to different equivalence classes is unique
(see [F, 5.3] for the compact group case). Also note that we do
not claim that irreducible \repns of compact groupoids are finite
dimensional. This is indeed false. I am indebted to Paul Muhly who
reminded me of the following simple counterexample: Let
$\G=[0,1]\times [0,1]$ with the Haar system
$\lambda^u=\delta_u\times \lambda$, where $\delta_u$ is the Dirac
point mass measure and $\lambda$ is the Lebesgue measure on
$[0,1]$. Then up to unitary equivalence) $G$ has only one
irreducible \repn $\pi$ with $\pi(C^*(\G,\lambda^u))=\mathcal
K(L^2[0,1])$, where the left hand side is the image of the
groupoid $C^*$-algebra of $\G$ (see [P] for details) and the right
hand side is the algebra of compact operators. This shows that
$\pi$ is not finite dimensional. On the other hand the bundle
giving $\pi$ is just the trivial bundle with one dimensional
fibers over $[0,1]$. The finite dimensionality of all irreducible
\repns is a very strong restriction on a groupoid $\G$. It is
believed that such a condition forces $\G$ to be a (measure
theoretic) bundle of transitive groupoids, each with a finite unit
space and compact isotropy groups [M].

Let $\pi\in\R$ and $\rho\in\hat\G$ be an irreducible sub\repn of
$\pi$. For $u\in X$, let $(\M_\rho)_u$ be the closed linear span
of all irreducible subspaces of $\Hu$ on which $\pi$ is equivalent
to $\rho$. Put $\M_\rho=\{(\M_\rho)_u\}_{u\in X}$. Then the
following is proved as [F, 5.3]

\begin{lem}
Let $\pi\in\R$ and $\rho,\rho_1,\rho_2\in\hat\G$ are irreducible sub\repns of $\pi$, then

$(i)$ $\M_{\rho}$ is invariant  under $\pi$ and if $\mathcal N$ is
any $\pi$-irreducible subbundle of $\M_\rho$, then $\pi^\mathcal
N\sim \rho$.

$(ii)$ If $\rho_1,\rho_2$ are not  unitary equivalent , then
$\M_{\rho_1}\perp\M_{\rho_2}$.\qed
\end{lem}

\begin{cor}
With above notation,
$$\HH\simeq\bigoplus_{\rho\in\hat\G} \M_\rho\simeq \bigoplus_{\rho\in\hat\G}
(\bigoplus_{\alpha\in\Lambda_\rho} \mathcal N_{\rho,\alpha}),$$
where
$$\pi^{\mathcal N_{\rho,\alpha}}\sim \rho\quad(\alpha\in\Lambda_\rho,\rho\in\hat\G).\qed$$
\end{cor}

\vspace{.3 cm} Note  that the second decomposition is not unique,
but as each $\rho\in\hat \G$ has a finite dimensional bundle, it
is trivial that the cardinality of $\mathcal N_{\rho,\alpha}$ is
independent of the decomposition and just depends on $\pi$ and
$\rho$. This is denoted by $mult(\rho,\pi)$ and is called the {\it
multiplicity} of $\rho$ in $\pi$. As in [F, 5.4] we have

\begin{lem}
With above notation, $mult(\rho,\pi)=dim Mor(\rho,\pi)$.\qed
\end{lem}

Note that when $\G$ is not compact , one can define the
multiplicity of irreducible sub\repns using the above equality
(which might result in an infinite cardinal) and Lemma 2.17 would
hold in locally compact case as well. Note that the decomposition
of Corollary 2.18 may fail in non compact case.

\section{Harmonic Analysis on compact groupoids}

In this section we turn into two very important results which are
of crucial importance in our duality theorem [A2], namely the
Peter-Weyl theorem and Schur's orthogonality relations. We start
with the definition of matrix elements.

Let $\pi\in\R$, The  mappings $x\mapsto
<\pi(x)\xi_{s(x)},\eta_{r(x)}>\quad (\xi,\eta\in\HH)$ are called
{\it    matrix elements} of $\pi$. This terminology is based on
the fact that if $\{e_u^i\}_i$ is a basis for $\Hu$, then
$\pi_{ij}(x)=<\pi(x)e_{s(x)}^j,e_{r(x)}^i>$ is the $(i,j)$-th
entry of the (possibly infinite) matrix of $\pi(x)$. We denote the
linear span of matrix elements of $\pi$  by $\EE$. By Lemma 2.3,
$\EE$ is a subspace of $C(\G)$. It is clear that $\EE$ depends
only on the unitary equivalence class of $\pi$. These vector
spaces are the building blocks of the Peter-Weyl theorem, so we
would like to have a closer look to their properties. First we
establish some notations and remind some facts.

Let's for a while go  back to the general case of a locally
compact groupoid $\G$. Each \repn $\pi$ of $\G$ could be
integrated to a \repn of the convolution algebra $C_c(\G)$ on
$L^2(\G^{(0)},\HH,\mu_\pi)$ via
$$<\pi(f)\xi,\eta>=\int_\G f(x)<\pi(x)\xi(s(x)),\eta(r(x))>d\nu_\pi(x)\quad (f\in C_c(\G),$$
where $\nu_\pi=\int_X \lambda_u d\mu_\pi(u)$. This  could also be
considered as a representation on the bundle $\HH$. Let's denote
the matrix element of $\pi$ at $\xi,\eta\in\HH$ by
$\pi_{\xi,\eta}$, namely
$\pi_{\xi,\eta}(x)=<\pi(x)\xi_{s(x)},\eta_{r(x)}>, \, x\in\G$. For
each function $\phi:\G\to\mathbb C$, let
$$L_x(\phi)(y)=\phi(x^{-1}y)\quad (y\in\G^{r(x)}),\,
R_x(\phi)(y)=\phi(yx)\quad(y\in\G_{r(x)}),$$ and as before, let
$\check \phi(y)=f(y^{-1}), \, y\in\G$. For each $\xi\in\HH$,
$x\in\G$, let $\pi(x)\xi$ denote the vector in $\HH$ whose fiber
at $u\in X$ is $\pi(x)\xi_u$, if $u=s(x)$, and $0$, otherwise.
Next lemma is valid for any locally compact groupoid. The proof is
straightforward and is omitted.

\begin{lem}
Let $\pi\in\R$, $\xi,\eta\in\HH$, $x\in \G$, and $f\in C_c( \G)$.  Then

$(i)\, L_x(\pi_{\xi,\eta})=\pi_{\xi,\pi(x)\eta}$ on $\G^{r(x)}$,

$(ii)\,  R_x(\pi_{\xi,\eta})=\pi_{\pi(x)\xi,\eta}$ on $\G_{r(x)}$,

$(iii)\, f*\pi_{\xi,\eta}=\pi_{\xi,\pi(\bar f)\eta}$ on $\G$,

$(iv)\, \pi_{\xi,\eta}*f=\pi_{\pi(\check f)\xi,\eta}$ on $\G$.

In particular $\EE$ is a two sided ideal of $C_c(\G)$, closed
under translations. \qed
\end{lem}

\vspace{.3 cm} There is a technical difficulty when one wants to
deal with $\EE$. Even if all Hilbert spaces $\Hu$ are finite
dimensional, there is no guarantee that $\EE$ is also finite
dimensional. Indeed, with above notation, if we write $\xi_{s(x)}
,\eta_{r(x)}$ as linear combinations of $e_{s(x)}^j$'s and
$e_{r(x)}^i$'s , the coefficients depend on $x$ in general, so
$\pi_{\xi,\eta}$ would not be a linear combination of
$\pi_{ij}$'s. To overcome this difficulty, we have to fix the
domain and range as follows. Take any $u,v\in X$ and restrict
$\pi_{\xi,\eta}$ to $\G_u^v$. As before let $\{e_u^i\}_i$ be a
Hamel basis for $\Hu$ (finite or infinite). decompose $\xi_u$ and
$\eta_v$ into finite linear combinations of basis elements. Then
$\pi_{\xi,\eta}$ is clearly a finite linear combination of the
corresponding $\pi_{i,j}$'s on $\G_u^v$. If  one uses a Schauder
basis for $\Hu$, the problem of convergence of the later
decomposition should be resolved. In any case, when $\Hu$ and
$\Hv$ are finite dimensional we get the following. For $u,v\in X$,
$\Euv$ consists of restrictions of elements of $\EE$ to $\G_u^v$.

\begin{prop}
Let $\pi\in\R$, $u,v\in X$, and $f$ be a  complex valued function
on $\G_u^v$. Assume that the Hilbert  spaces $\Hu$ and $\Hv$ are
finite dimensional. Then the following are equivalent.

$(i) f\in\Euv$,

$(ii)$ There is $A\in\mathcal B(\Hv,\Hu)$ such that $f=Tr(A\pi(.))$,

$(iii)$ There are bases $\{e_u^j\}_{1\leq j\leq d_u^\pi}$  and
$\{e_v^i\}_{1\leq i\leq d_v^\pi}$ of $\Hu$ and $\Hv$ such that $f$
is a linear combination of the matrix elements
$\pi_{u,v}^{ij}=<\pi(.)e_u^j,e_v^i>$,

$(iv)$ There are $\xi^j,\eta^i\in\HH\quad (1\leq j\leq
d_u^\pi,\,1\leq i\leq d_v^\pi)$  such that
$$f=\sum_{i=1}^{d_v^\pi}\sum_{j=1}^{d_u^\pi}
<\pi(.)\xi_u^j,\eta_v^i>.$$
\end{prop}
{\bf Proof} The bounded linear functionals on $\mathcal
B(\Hu,\Hv)$  are exactly the maps $B\mapsto Tr(AB)$, where
$A\in\mathcal B(\Hv,\Hu)$, so $(i)$ and $(ii)$ are equivalent.
Equivalence of $(ii)$ and $(iii)$ is an exercise in elementary
linear algebra. The equivalence of $(iii)$ and $(iv)$ is
trivial.\qed

\begin{cor}
If $\Hu$  and $\Hv$ are finite dimensional, then so is $\Euv$.
Moreover $dim(\Euv)\leq (dim\Hu).(dim\Hv)$.\qed
\end{cor}
Next we prove a technical (but easy) lemma which is of crucial
importance in proving the Peter-Weyl theorem for groupoids.

\begin{lem} Let $\pi_1,\pi_2\in\R$ and $A:\HH\to\HH$ is any bundle of
bounded linear operators, put

$$\tilde A_u\xi_u =\int \pi_2(x^{-1})A_{r(x)} \pi_1(x)
d\lambda_{u}(x),\quad (u\in X,\xi\in\HH)$$ then $\tilde A\in
Mor(\pi_1,\pi_2)$.
\end{lem}
{\bf Proof} Given $x\in \G$,
\begin{align*}
\tilde A_{r(x)}\pi_1(x) &=\int \pi_2(y^{-1})A_{r(y)} \pi_1(yx) d\lambda_{r(x)}(y) \\
&=\int \pi_2(xy^{-1})A_{r(y)} \pi_1(y) d\lambda_{s(x)}(y)=\pi_2(x) \tilde A_{s(x)}.\qed
\end{align*}

\begin{lem} If $\pi_1,\dots,\pi_n\in\R$ and $\pi=\pi_1\oplus\dots\oplus\pi_n$,
then $\Euv=\sum_{i=1}^{n} \mathcal E_{u,v}^{\pi_i}$ (not a direct sum).
\end{lem}
{\bf Proof} Proof is a straightforward calculation.\qed

\begin{thm} {\bf (Schur's orthogonality relations)} Let $u,v\in X$ and $\pi,\pi^{'}\in\R$.
Consider $\Euv$ and $\mathcal E_{u,v}^{\pi^{'}}$ as subspaces of $L^2(\G_u^v,\lambda_u^v)$.

$(i)$ If $\pi,\pi^{'}\in\R$ and $\pi\nsim \pi^{'}$, then $\Euv\perp\mathcal E_{u,v}^{\pi^{'}}$,

$(ii)$ If $\lambda_u(\G_u^v)\neq 0$, then $dim(\Euv)=d_u^\pi d_v^\pi$ and
$\{\sqrt{d_u^\pi} \pi_{u,v}^{ij}: 1\leq i\leq d_v^\pi, 1\leq j\leq d_u^\pi\}$ is
an orthonormal basis for $\Euv$.
\end{thm}
{\bf Proof} Let $\eta\in\HH,\eta^{'}\in\mathcal H_{\pi^{'}}$. In
Lemma 2.22, put $A_u\xi_u=<\xi_u,\eta_u>\eta_u^{'}$, then  for
each $u\in X$, $\xi\in\HH,\xi^{'}\in\mathcal H_{\pi^{'}}$,
\begin{align*}
<\tilde A_u\xi_u,\xi_u^{'}>&=\int <\pi^{'}(x^{-1})A_{r(x)}\pi(x)\xi_u,\xi_u^{'}>d\lambda_u(x)\\
&=\int <A_{r(x)}\pi(x)\xi_u,\pi^{'}(x)\xi_u^{'}>d\lambda_u(x)\\
&=\int <\pi(x)\xi_{s(x)},\eta_{r(x)}><\eta_{r(x)}^{'},\pi^{'}(x)\xi_{s(x)}^{'}>d\lambda_u(x)\\
&=\int \pi_{\xi,\eta}(x)\overline{\pi_{\xi^{'},\eta^{'}}(x)} d\lambda_u(x).
\end{align*}
Now by Schur's lemma applied to $\tilde A\in Mor(\pi,\pi^{'})$,
if $\pi\nsim\pi^{'}$, then  $\tilde A=0$, and so $(i)$ follows
from above equalities. Again by Schur's lemma, if $\pi=\pi^{'}$,
$\tilde A=cI$, so if we take $\xi=e^i,\eta=e^j,\xi^{'}=e^{i^{'}},$
and $\eta^{'}=e^{j^{'}}$, then
$$\int \pi_{ij}(x)\overline{\pi_{i^{'}j^{'}}(x)} d\lambda_u(x)
=c<e_u^i,e_u^{i^{'}}>=c\delta_{ii^{'}}.$$
Also
$$cd_u^\pi=Tr(\tilde A_u)=\int Tr(A_{r(x)}\pi(x^{-1}x)d\lambda_u(x)
=\int Tr(A_{r(x)})d\lambda_u(x),$$
but for each $v\in X$,
$$Tr(A_v)=\sum_{i=1}^{d_v^\pi} <A_ve_v^i,e_v^i>=\sum_i
<e_v^i,e_v^j><e_v^{j^{'}},e_v^{i}>=\delta_{jj^{'}},$$
so $cd_u^\pi=\delta_{jj^{'}}\lambda_u(\G_u)=\delta_{jj^{'}}.$ Therefore
$$<\pi_{ij},\pi_{i^{'}j^{'}}>=\frac{\delta_{ii^{'}}\delta_{jj^{'}}}{d_u^\pi},$$
as  functions on $\G_u$. Now we can easily redo all the above
calculations with all integrals, started with the one in the
definition of $\tilde A$, taken over $\G_u^v$ (instead of $\G_u$).
The only difference would be the value of the constant $c$, as
this time we get $cd_u^\pi=\delta_{jj^{'}}\lambda_u(\G_u^v)$. Now
if $\lambda_u(\G_u^v)=0$, then clearly $\Euv=\{0\}$, otherwise,
$$<\pi_{u,v}^{ij},\pi_{u,v}^{i^{'}j^{'}}>
=\frac{\delta_{ii^{'}}\delta_{jj^{'}}}{d_u^\pi\lambda_u(\G_u^v)},$$
and we get $(ii)$.\qed

\vspace{.3 cm} Next,  for $u,v\in X$, consider the left and right
regular representations $L$ and $R$ acting on
$L^2(\G_u^v,\lambda_u^v)$.

\begin{lem}
Let $\pi\in\hat\G$ and put
$$\mathcal R_{u}^i=span\{\pi_{u,v}^{i1},\dots,\pi_{u,v}^{id_u^\pi}\}
\quad(1\leq i\leq d_v^\pi),$$
and
$$\mathcal C_{v}^j=span\{\pi_{u,v}^{1j},\dots,\pi_{u,v}^{d_v^\pi j}\}
\quad(1\leq j\leq d_u^\pi).$$ Then  each $\mathcal R_{u,v}^i$
($\mathcal C_{u,v}^j)$ is invariant under right (left) regular
\repn and $R^{\mathcal R_{u,v}^i}\sim \pi$ ($L^{\mathcal
C_{u,v}^j}\sim\bar\pi$).
\end{lem}
{\bf Proof}
Let $\{e_u^i\}$ be an orthonormal basis for $\Hu$. Then for each $x\in \G_u^v$,
$$\pi(x)(\sum_j c_je_u^j)=\sum_j c_j\sum_k <\pi(x)e_u^j,e_v^k>e_v^k
=\sum_{j,k} c_j\pi_{u,v}^{kj}(x)e_v^k,$$ now  if $y\in \G_v^w$,
then one can calculate $\pi(yx)(\sum_j c_j e_u^j)$ in two
different ways, by substituting $x$ by $yx$ in above formula or by
applying $\pi(y)$ to the both side of that formula. This gives us
two different decomposition into linear combinations of the basis
elements $e_w^i$ and if we put the coefficients equal we get
$$\pi_{u,w}^{ij}(yx)=\sum_k \pi_{v,w}^{ik}(y)\pi_{u,v}^{kj}(x), $$
therefore
$$R_x(\sum_j c_j\pi_{u,w}^{ij})=\sum_{j,k} c_j\pi_{u,v}^{kj}(x)\pi_{v,w}^{ik}.$$
This shows that if $U_u^i:\Hu\to\mathcal R_u^i$ is defined by
$$U_u^i(\sum_j c_je_u^j)=\sum_j c_j\pi_{u,v}^{ij},$$
then $U^i$ is clearly a bundle of unitaries and
$U_v^i\pi(x)=R_xU_u^i\quad(x\in\G_u^v)$, that is $U\in
Mor(\pi,R^{\mathcal R_{u,v}^i})$ and the first statement is
proved. The proof of the other statement is similar.\qed

\vspace{.3 cm} Let us put
$$\mathcal E_{u,v}=span\big(\cup_{\pi\in\hat\G}\Euv\big)$$
 and
$$\mathcal E=span\big(\cup_{u,v\in X} \mathcal E_{u,v}\big).$$

\begin{prop}
$\mathcal E_{u,v}$ and $\mathcal E$ are algebras.
\end{prop}
{\bf Proof} Let $\pi,\pi^{'}\in\hat\G$ and choose bases
$\{e_u^i\}$ and  $\{f_u^j\}$ for $\Hu$ and $\mathcal
H_u^{\pi^{'}}$, respectively. For $x\in\G_u^v$, define
$\pi\otimes\pi^{'}(x):\Hu\bigotimes\mathcal H_u^{\pi^{'}}\to
\Hv\bigotimes\mathcal H_v^{\pi^{'}}$ by
$$\pi\otimes\pi^{'}(x)(\xi_u\otimes\xi_u^{'})
=\pi(x)(\xi_u)\otimes \pi^{'}(x)(\xi_u^{'})
\quad(\xi_u\in\Hu, \xi_u^{'}\in\mathcal H_u^{\pi^{'}}).$$
Then for all indices $i,j,k$, and $l$,
$$
<\pi\otimes\pi^{'}(x)(e_u^j\otimes f_u^l,(e_v^i\otimes f_v^k)> =
<\pi(x)e_u^j,e_v^i>.<\pi^{'}(x)f_u^l,f_v^k>,$$ that  is
$(\pi\otimes\pi^{'})_{u,v}^{ikjl}=\pi_{u,v}^{ij}{\pi^{'}}_{u,v}^{kl}$,
which proves the first assertion. The proof for $\mathcal E$ is
similar.\qed

\vspace{.3 cm} Next we prove a result which is the main ingredient
of both Gelfand-Raikov and Peter-Weyl theorem for compact
groupoids. The proof closely follows the original proof of Peter
and Weyl [F, 5.11].

\begin{lem}
$\mathcal E$ is dense in $C(\G)$ .
\end{lem}
{\bf Proof}
Consider the following operators
$$T_\psi(f)(x)=\psi *f(x)=\int \psi(xy^{-1})f_{s(x)}(y)d \lambda_{s(x)}(y)
\quad (x\in \G,\psi\in C(\G), f\in L^2(\G)),$$ where $L^2(\G)$  is
the bundle whose fiber at $u\in X$ is $L^2(\G_u,\lambda_u)$ ( and
so each $f\in L^2(\G)$ is of the form $f=\{f_u\}$, with $f_u\in
L^2(\G_u,\lambda_u)$ and
$\|f\|_2=sup_u\|f_u\|_{L^2(\G_u,\lambda_u)}<\infty$). In other
words, $T_\psi(f)(x)=L_{s(x)}(\psi)(f_{s(x)})$, where $L$ is the
left regular \repn of $C(\G)$ on $L^2(\G)$. It is easy to see that
$T_\psi^{*}=T_{\check \psi}$, so if $\psi$ is symmetric, $T_\psi$
is self-adjoint. Also it is routine to check that $T_\psi(f)\in
C(\G)$ and $\|T_\psi(f)\|_\infty \leq \|\psi\|_\infty\|f\|_2$.
Also if $\ell_x$ is the operator of left translation by $x\in \G$
on $C(\G)$, then
$$\|\ell_x(T_\psi(f))-T_\psi(f)\|_\infty =\|(\ell_x(\psi)-\psi)*f\|_\infty
\leq \|\ell_x(\psi)-\psi\|_\infty\|f\|_2,$$ so if $B$ is a bounded
set in $L^2(\G)$, $T_\psi(B)$ is uniformly bounded and
equicontinuous in $C(\G)$, and so it follows from Arzela-Ascoli
theorem that $T_\psi$ is a compact operator on $L^2(\G)$. One can
consider $T_\psi$ as an operator on the fiber
$L^2(\G_u,\lambda_u)$ and by the same argument, it would be a
compact operator. Now by an standard theorem of operator theory
[F, 1.52], $L^2(\G_u,\lambda_u)$ has an orthonormal basis
consisting of eigenvectors of $T_\psi$. Given an eigenvalue $z$ of
$T_\psi$, the corresponding eigenspace $\mathcal M_z$ consists of
those $g\in L^2(\G_u,\lambda_u)$ with $\psi*g=zg$. Now if this
holds for $g$ and $r_x$ is the right translation by $x$, then
$\psi*r_x(g) =r_x(\psi*g)=zr_x(g)$, that is the eigenspace of $z$
is right translation invariant, so it is invariant under the right
regular \repn $R$. $\mathcal M_z$ is clearly finite dimensional,
if $z\neq 0$. Choose an orthonormal basis $e_u^1,\dots ,e_u^{n_u}$
for $\mathcal M_z$ and let $R_u^{jk}$ be the corresponding
coefficient functions of the right regular \repn on $\mathcal
M_z$. Then for $u,v\in X$, and each $x\in \G_u^v,y\in\G_v$ and
$1\leq k\leq n$, $e_u^k(yx)=\sum_{i=1}^{n_v}  R_u^{ik}(x)
e_v^i(y)$, and so $e_u^k=\sum_i e_v^i(v)R_u^{ik}$, which shows
that $\mathcal M_z\subseteq \cup_{v\in X} \mathcal
E_{u,v}^R\subseteq \mathcal E$. Next if $g\in L^2(\G_u,\lambda_u)$
, then $g=\sum_z g_z$ with $g_z\in \mathcal M_z$, and the series
converges in $L^2(\G_u,\lambda_u)$. By boundedness of
$T_\psi:L^2(\G_u,\lambda_u)\to C(\G)$, $T_\psi(f)=\sum_z zg_z$,
and the series is converging in the uniform norm of $C(\G)$. By
what we have just shown, any finite sub summation of this series
is in $\mathcal E$, and so $\mathcal E\cap Im(T_\psi)$ is
uniformly dense in $Im(T_\psi)$, but the union of $Im(T_\psi)$'s
when $\psi$ ranges over an approximate identity of $C(\G)$
consisting of symmetric functions, is dense in $C(\G)$, so
$\mathcal E$ is dense in $C(\G)$.\qed

\begin{cor}
For each $u,v\in X$, $\mathcal E_{u,v}$ is dense in $C(\G_u ^v)$.\qed
\end{cor}

\vspace{.3 cm} The following result is proved as in [HR,1.8.4].
\begin{lem}
Each Hausdorff locally compact groupoid is completely Hausdorff (as a topological space).
\end{lem}

\vspace{.3 cm} Now we are ready to prove two of the main results
of this paper.

\begin{thm} {\bf (Gelfand-Raikov Theorem)}
If $\G$ is a Hausdorff compact groupoid, then $\hat \G$ separates the points of $\G$.
\end{thm}
{\bf Proof} Let $x,y\in\G$ and $x\neq y$. By the above  lemma,
there is a function $f\in C(\G)$ such that $f(x)\neq f(y)$. By
Lemma 3.9, we might assume that $f\in\mathcal E$. Then $f$ is a
finite linear combination of the coefficient functions of some
elements of $\hat \G$. Therefore for at least one of these, say
$\pi\in\hat\G$, we have
$<\pi(x)\xi_{s(x)},\eta_{r(x)}>\neq<\pi(y)\xi_{s(y)},\eta_{r(y)}>$,
for some unit vectors $\xi,\eta\in\HH$. Hence we have
$\pi(x)\neq\pi(y)$, as required.\qed

\vspace{.3 cm} Summing up the results of Theorem 3.6, Lemma 3.7,
and Corollary 3.10, we have

\begin{thm} {\bf (Peter-Weyl Theorem)}
Let $\G$ be a compact groupoid, then for each $u,v\in X$,  $\Euv$
is dense in $C(\G_u^v)$ and
$$L^2(\G_u^v,\lambda_u^v)=\bigoplus_{\pi\in\hat\G} \Euv,$$
and $$\{\sqrt{d_u^\pi}\pi_{u,v}^{ij}: \pi\in\hat\G,\, 1\leq i\leq
d_v^\pi, \,1\leq j\leq d_u^\pi\},$$ is an orthonormal basis for
$L^2(\G_u^v,\lambda_u^v)$. Each  $\pi\in\hat \G$ occurs in the
right and left regular \repn of $\G$ over
$L^2(\G_u^v,\lambda_u^v)$ with multiplicity $d_u^\pi$. \qed
\end{thm}


\end{document}